\documentclass[]{article}
\usepackage[utf8]{inputenc}
\usepackage{natbib}
\usepackage{tabularx, booktabs}
\usepackage{multirow}
\usepackage{url}
\usepackage[colorinlistoftodos]{todonotes}
\usepackage{listings}
\usepackage{xcolor}
\usepackage{nameref}
\usepackage{mathtools}
\usepackage{comment}
\definecolor{codegreen}{rgb}{0,0.6,0}
\definecolor{codegray}{rgb}{0.5,0.5,0.5}
\definecolor{codepurple}{rgb}{0.58,0,0.82}
\definecolor{backcolour}{rgb}{0.95,0.95,0.92}
\lstdefinestyle{mystyle}{
    backgroundcolor=\color{backcolour},   
    commentstyle=\color{codegreen},
    keywordstyle=\color{magenta},
    numberstyle=\tiny\color{codegray},
    stringstyle=\color{codepurple},
    basicstyle=\ttfamily\footnotesize,
    breakatwhitespace=false,         
    breaklines=true,                 
    captionpos=b,                    
    keepspaces=true,                 
    numbers=left,                    
    numbersep=5pt,                  
    showspaces=false,                
    showstringspaces=false,
    showtabs=false,                  
    tabsize=2
}
\title{PathWise: a flexible, open-source library for the Resource Constrained Shortest Path}
\author{Matteo Salani and  
Saverio Basso\\Dalle Molle Institute for Artificial Intelligence (IDSIA)\\USI-SUPSI, Lugano, Switzerland}

\begin{document}

\maketitle

\begin{abstract}
In this paper, we consider a fundamental and hard combinatorial problem: the Resource Constrained Shortest Path Problem (RCSPP).
We describe the implementation of a flexible, open-source library for the solution of the RCSPP, called \textit{PathWise}, capable of tackling several variants of the problem.
We designed PathWise with the final user in mind, developing easy-to-use interfaces without compromising performance. 
We provide computational experiments on three classes of instances of the RCSPP, namely RCSPP on cyclic networks, RCSPP on large acyclic networks, and RCSPP on ad-hoc cyclic networks. We show that PathWise is packed off-the-shelf with algorithms capable of tackling all classes.
This paper represents the first step along the journey of devising and implementing a comprehensive open-source library for a large variety of RCSPPs. The current version of the library carries exact algorithms for the RCSPP but new algorithms, both heuristic and exact, will be added thanks to the flexible design. 
We also foresee PathWise becoming a platform ready for data-driven and process-driven methodologies for these types of problems. 
\end{abstract}

\section{Introduction}

Shortest path problems are everywhere. Their application goes far beyond the computation of the shortest route in our favorite road navigation system.
They are at the cornerstone of the most challenging and practically relevant combinatorial problems in transportation, telecommunication, and scheduling. Shortest paths also appear as subproblems in routing and workforce planning problems solved by column generation \citep{Irnich2005}.

While it is commonly believed that computing the shortest path is an easy problem, this is only true in specific cases: when the shortest path is unconstrained and the underlying graph does not possess negative cost cycles. In the large majority of relevant applications, either the shortest path must respect some additional constraints, commonly referred to as resource constraints, or the network exhibits negative cost cycles or both cases happen. This makes the Resource Constrained Shortest Path Problem (RCSPP) a hard combinatorial problem. 

The RCSPP is defined on a graph $G(N,A)$, which we assume directed, composed of a node set $N$ and an arc set $A$. The problem asks to find a minimum cost elementary path, i.e. a finite sequence of consecutive arcs in which every node $n \in N$ appears at most once, from a source node $s \in N$ to a destination node $d \in N$. The cost is accumulated when traversing arcs along the path. We remark that no assumptions are taken on the cost of the arcs and the graph may possess negative cost cycles. 

The RCSPP problem has one or more resources that are consumed while traversing arcs. For example, elapsed time, transported load, etc., and their availability is constrained.
In some problems, resources can be more complex and their availability can be renewed at some cost, for example, the energy stored in a battery of an electric vehicle can be restored by spending time at a recharging station.

In this paper, we describe the design and implementation of \textit{PathWise} a flexible open-source library for the solution of the RCSPP.
PathWise has been designed with both beginner and expert users in mind: beginners can solve a variety of standard RCSPPs with an off-the-shelf implementation of state-of-the-art algorithms by calling PathWise from their favorite programming language, and experienced users are capable of developing their algorithmic components of PathWise while taking advantage of the framework thanks to clear interfaces and well-defined hook points.  

PathWise has another long term objective in mind: to be a playground for the integration of Machine Learning techniques with heuristic and exact algorithms for several classes of RCSPP-like problems. 

\section{Literature}
A great number of real-life applications can be modeled with resource constrained shortest paths formulations. Some examples include vehicle and crew scheduling problems \citep{Desaulniers1998,Haase2001}, rostering \citep{Gamache1999},  military aircraft management systems \citep{Zabarankin2002}, railroad management \citep{Halpern1974}, telecommunication network design \citep{Cabral2007}, green vehicle routing problems \citep{ERDOGAN2012} and many others. 

In the following, we refer to the works of \citep{Irnich2005} and \citep{Pugliese2013} for a more detailed survey of the literature concerning the main contributions for the RCSPP and its variants. A more general review of Shortest Path Problems can be found in \citep{Madkour2017}.

Some classical approaches for solving RCSPP include pre-processing procedures and Branch \& Bound \citep{Christofides1989} and Lagrangian relaxations and enumeration of near-shortest paths \citep{Carlyle2008}. In particular, \emph{Dynamic Programming labeling algorithms} \citep{Mehlhorn2000, Dumitrescu2003,  Desrochers1988} are among the most successful exact methods. They are based on labels that describe the state of the algorithm, that is, they encode partial path information from the source node to another node of the network. Propagation and management of only the pareto-optimal labels associated at every node, allows these methods to solve very efficiently well sized instances but, since performance is strongly dependent on the number of generated labels, they often struggle when facing large scale networks.

Notable improvements can be obtained by including bidirectional search \citep{RighiniSalani2006}, by propagating labels from both the source to the destination and backward, from the destination to the source, and by joining partial paths. The idea has been further improved by dynamically balancing forward and backward label pools in \citep{Tilk2017}.

Furthermore, recent techniques exploit bucket based approaches \citep{PecinPessoa2017, Sadykov2020}. In particular, in \citep{Sadykov2020}, the authors propose a variant of a dynamic programming algorithm based on a bucket graph: in this version, labels are collected and extended in buckets. This has the main effect of decreasing the number of comparisons for dominance, resulting in significant improvements in running time.

Many heuristics can also be found in the literature. For example, in \citep{Desaulniers2008} the authors relax dynamic programming dominance rules by focusing only on a subset of the available resources. Other attempts, include limiting the number of labels that can be stored at each node, like in \citep{Fukasawa2006}. 

While resource constraints are common characteristics of different variants of the RCSPP, a distinctive feature is the presence of negative cost cycles in the underlying network, generally modeling problems where nodes have associated prizes. Most relevant literature contributions can be indeed classified as methods specialized in acyclic or cyclic networks.

\paragraph*{Methods for acyclic networks.} During the last decade, dynamic programming techniques for acyclic graphs have been proposed. For example, \emph{Pulse} \citep{Lozano2013} and  \emph{bidirectional Pulse} \citep{Lozano2020} algorithms explore the graph with a depth-first search strategy,  along with pruning, parallelism options and opportunities to redirect the search towards more promising solutions. In \citep{RCBDA} the authors combine instead, elements bidirectional A* and bidirectional dynamic programming, to obtain an exact algorithm. These techniques resulted in  performing on relatively large sized instances. Very recently, an enhanced, biased A* algorithm has been proposed in \citep{EBBA}, along with several heuristics, improving the speed and quality of bound computation, and pruning procedures, the balance forward and backward extension steps, and the efficiency of dominance and join methods. Overall, the algorithm quickly obtains solutions for shortest path problems with a single resource over large scale networks.

\paragraph*{Methods for cyclic networks.} The more complex variant of the problem however, the \emph{Resource Constrained Elementary Shortest Path Problem (RCESPP)}, arises when facing cyclic graphs, and requires finding elementary paths, i.e. with no repeating nodes. It is NP-hard in the strong sense and has been the focus of numerous studies.

 Some applications include Orienteering problems \citep{Golden1987} \citep{Gunawan2016}, that search for a path between two nodes that maximizes prize collection, and generalizations of the Travelling Salesman problems \citep{Laporte1990, Feillet2005}, that target visiting a subset of the nodes, with simultaneous optimization of profits and travel costs.
Furthermore, related problems such as Vehicle Routing Problems and Crew Scheduling formulations are usually solved with Column Generation \citep{Desaulniers2005}. In this setting, pricing problems generally correspond to a RCESPP with prize collection, whilst the restricted master problem can be modeled as a set partitioning that selects the most promising paths.
An effective Branch-and-Cut algorithm has been proposed in \citep{spprclib}.
However, most of the research exploits relaxations to overcome the complexity of RCESPP. Among the top techniques, in \citep{RighiniSalani2008} the authors propose to progressively lift the state space, by iteratively populating a set of nodes necessary to compute an elementary solution. In this direction, \emph{ng-path} relaxations have been proposed in \citep{Baldacci2011}: they define, for each node, a set of neighbor nodes that must not contain cycles. Efficient implementations of these relaxation methods have been proposed in \citep{Martinelli2014}. 
 Very recently, variants of ng formulations based on managing sets of neighbor arcs instead of nodes have been introduced in \citep{Bulhoes2018} and \citep{Costa2021}. They proved effective in reducing the number of non-dominated labels and therefore computing times. 
Other relevant approaches such as \citep{Irnich2006} and \citep{Desaulniers2008}, rely respectively on forbidding cycles of small length and on relaxing elementarity requirements for selected nodes. 
Finally, heuristics for RCESPP have been studied as well, like in \citep{HOMBERGER2005}. 

Regarding Branch \& Price, some methods to also improve the restricted master problem performance have been proposed. In \citep{Jepsen2008}, the authors introduce in the formulation new subset-row inequalities, that is additional Chvatal rank-1 cuts, and, for each of them, consider an additional resource in the labeling algorithm for RCESPP. This has the side effect of making dominance more difficult, albeit better lower bounds can be obtained and, in some cases, optimality at the root node can be proven. Limited memory techniques can be used to weaken these cuts to mitigate negative effects on pricing \citep{PecinPessoa2017}.

\paragraph*{Parallel computing and algorithm collections.} From the technological point of view, standard parallel implementations of dynamic programming algorithms are possible. Recent solutions include \citep{Lu2020}, in which the authors propose a framework for parallelization of labeling algorithms through GPU, reporting impressive speed-ups.

We also report that some collection of algorithms for the RCSPP and related problems can be found in the literature but, at this time, none of them includes recent state of the art techniques. For example, the authors of \citep{Sanchez2020} present a Python library that employs bidirectional labeling and metaheuristics for the RCSPP with multiple non-decreasing resources. Collections of algorithms for different but related problems, include \citep{Montagne2020}, a package for VRP that exploits column generation, albeit with no optimality guarantees.
A generic Branch-Cut-and-Price exact solver for VRP and related problems is instead proposed in \citep{Pessoa2020}, also exploiting other state-of-the-art techniques such as ng-paths, variants of subset row and capacity cuts, path enumeration and bucket graphs. It performs comparably against ad-hoc VRP implementations.

\section{PathWise} \label{PW}

The first release of PathWise ships state-of-the-art algorithms based on dynamic programming, along with tailored contributions on relaxation schemes and bidirectional dynamic programming. We collected such algorithms in a single framework with an efficient implementation.
We report hereafter the basics of these algorithms, describe the conceptual design, and provide some use cases of the library.

\subsection{Algorithms} \label{Algorithms}

In dynamic programming, a state associated with node $i \in N$ represents a partial path from the source node $s$ to the node $i$. Different states can be associated with the same node and they correspond to different partial paths. The algorithm iteratively extends states until no further extensions are possible. Among all feasible states reaching the destination node $d$ the one with minimal cost represents the optimal solution to the RCSPP. Bi-directional dynamic programming simultaneously considers {\it forward} partial paths from the source node $s$ and {\it backward} partial paths reaching the destination $d$.

Each state is encoded in a label, in bi-directional dynamic programming called {\it forward} and {\it backward} labels. A forward label associated with node $i \in N$ is a tuple:
\begin{equation}
    l^f_i = (i, c_i, S, R)
\end{equation}
where $i$ is the last node visited in the partial path, $c_i$ is the accumulated cost, $S$ is a binary vector that keeps track of the visited nodes in the partial path and $R$ is the so-called resource vector that accounts for the consumption of each resource. 
Similarly, a backward label associated with node $i \in N$, corresponds to paths from node $i$ to destination node $d$. To control the number of labels, dominance tests are performed. Dominated labels can be safely discarded as they will not lead to an optimal solution.

The extension of a forward label corresponds to appending an additional arc $(i, j)$ to a path from $s$ to $i$, obtaining a path from $s$ to $j$, while the extension of a backward label corresponds to pre-pending an additional arc $(j, i)$ to a path from $i$ to $d$, obtaining a path from $j$ to $d$. Labels are stored in convenient data structures, referred to as {\it label pools}.

When a label $l_i = (i,c_i, S, R)$ is extended to a node $j$, a new label $l_j = (j, c_j, S', R')$ is generated by setting the corresponding element $S'_j$ to $1$ and updating the resource consumption vector $R'$ accordingly. 

Bi-directional dynamic programming largely reduces the number of generated labels by selecting a monotone resource, called {\it critical resource}, and extending labels for which its consumption is less than a given threshold $T$. 

Once the extension process is completed, forward and backward labels are joined to produce complete paths from node $s$ to node $d$. The join condition ensures that the final path contains no cycles nor violates resource constraints. Among the feasible paths that are generated by the join operation, we compute the optimal path. Its existence is guaranteed by the domination criteria.

When the underlying network is acyclic (i.e., the cost matrix does not possess negative cost cycles), all partial paths with cycles are suboptimal and the binary vector $S$ can be safely removed from the state vector. In case of cyclic networks, instead, the path elementarity does not come for free from cost minimization but must be specifically enforced using the binary vector $S$. Therefore, in a basic implementation of a dynamic programming algorithm for the RCESPP, the size of the state space grows exponentially with the size of the network.

Recent studies have worked on improving dynamic programming algorithms in multiple directions:  in this first implementation of the library we mainly focused on relaxation approaches and bidirectional techniques.

\paragraph*{Relaxation schemes} Decremental state space relaxation (DSSR, \citep{RighiniSalani2008}) aims at reducing the number of states to be explored by dynamic programming. 
The basic idea is that only a subset of the nodes of the network are relevant to compute the optimal solution without cycles. Therefore, the binary vector $S$ should be restricted only to those nodes. As this set is unknown, the algorithm starts with an empty set, relaxing a large portion of the state space, and iteratively adds nodes until a solution without cycles is found. Attempts to initialize the set $S$ have been explored in \citep{RighiniSalani2009}.

An alternative and very effective relaxation scheme has been proposed by \citep{Baldacci2011}. The relaxation, called {\it ng-path} relaxation, consists of defining, for each node of the graph, a subset of nodes on which path elementarity is enforced. The idea shares some similarities with that of the DSSR algorithm, but here the subsets at different nodes are independent. More recently, \citep{Martinelli2014} proposed to hybridize the ng-path relaxation with DSSR. In practice, the neigbourhood of each node on which elementarity is enforced is iteratively enlarged according to the optimal result of the relaxation until the optimal solution is an elementary shortest path.
In PathWise we implement different techniques that guarantee complete or partial elementarity. For the former, we implemented DSSR based methods:
\begin{itemize}
\item {\tt DSSR:} exploits the original DSSR scheme, here all nodes that are visited more than once in the optimal solution are added to the set $S$.
\item {\tt DSSRC:} stores a dedicated set for elementarity at each node and uses a DSSR based approach to iteratively forbid partial paths contaning cycles. Namely, when a cycle is found, nodes that appear more than once are inserted in the sets of nodes appearing in the loop only. For instance, let {\tt [0 5 6 8 6 9 5 3]} be the list of visited nodes in the optimal path of an intermediate step of the relaxation, {\tt DSSR} would add nodes {\tt 5} and {\tt 6} to the set $S$ of all nodes, forbidding any cycle over these two nodes. In contrast, {\tt DSSRC} would insert node {\tt 5} in the sets of nodes {\tt 5, 6, 8, 9} and node {\tt 6} to the sets of {\tt 6, 8} only. That is, loops over {\tt 5} and {\tt 6} might still happen in subsequent iterations, but not over the updated sets.

\end{itemize}
\noindent
For the latter, we implemented NG based procedures:
\begin{itemize}
\item {\tt NG:} uses the original NG-path approach to require partial elementarity. That is, path elementarity is enforced, at each node, for a set of neighbors only. 
\item {\tt NGC:} exploits a {\tt DSSRC} like iterative approach, relaxing the state space of the original {\tt NG} method. The algorithm starts with empty sets, and at each iteration, repeated visits to nodes that do not satisfy the original {\tt NG} scheme are forbidden. In this case, only the sets of nodes appearing in a loop are updated. 
\end{itemize}
We remark that both NG based techniques do not guarantee an elementary solution. However, at any time, hybridizations are possible:
\begin{itemize}
\item {\tt NG-DSSRC:} employs the {\tt NG} relaxation. If the solution presents cycles, enforces {\tt DSSRC} rounds to achieve an elementary path. 
\item {\tt NGC-DSSRC:} utilizes the {\tt NGC} relaxation first, followed by possible {\tt DSSRC} iterations.  
\end{itemize}

In PathWise, we propose an unified implementation to manage different relaxation schemes. In particular, we define for each node $i \in N$ a bit-mask $B_i = \{b^0, b^1, \dots b^N\}$ where $b^j \in \{0,1\}$. If $b^j = 1$ in the bit-mask $B_i$ of node $i$, it means that node $j$ belongs to the neighbourhood of node $i$, thus elementarity should be enforced.

When a label $l_i = (i,c_i, S, R)$ is extended to a node $j$, the set $S'_j$ is first initialized with $S'_j = S \& B_j$ (bit-wise \textit{And} operation) and then the $j-th$ element is set to $1$. In other words, the status ``forgets'' the previous visits to the nodes not beloning to the bit-mask $B_j$.

In order to implement DSSR schemes, the bit-masks of all nodes are equal, while for NG based schemes, the bit-masks are initialized with node-specific neighbourhoods.

\paragraph*{Bidirectional dynamic programming} 
Recently \citep{Tilk2017} proposed a technique, called \textit{dynamic half-way point}, that computes and adjusts the threshold $T$ during the label extension procedure in order to keep the dimensions of the forward and backward sets of labels as balanced as possible. 
In our implementation, we exploit instead the iterative nature of our relaxation algorithms.

Earlier iterations of any DSSR or NGC configurations are fast to solve and an uneven balance of forward and backward labels does not affect performance in a significant way. Anyway, after each iteration, the ratio between the overall number of forward ($N_F$) and backward ($N_B$) generated labels provides a good indication to compute a new half-way point ($hwp$), by taking into account the upper bound of the critical resource ($U_c$) as follows: 

\begin{equation}{
    hwp = 
    \begin{cases}
        hwp + 5\% \cdot U_c, & \text{if  }  \frac{N_B-N_F}{N_F} > 20\% \\
        hwp - 5\% \cdot U_c, &  \text{if } \frac{N_F-N_B}{N_B} > 20\%\\
        hwp, &\text{otherwise}\\
    \end{cases}
    }
\end{equation}

That is, if the number of generated forward (resp. backward) labels is more than 20\% the number of backward (resp. forward) labels we decrease (resp. increase) the half-way point accordingly to reduce the unbalance in subsequent iterations.
We call this technique \textit{Semi-dynamic half way point}. A comparable update strategy has been proposed in \citep{Sadykov2020}, for Vehicle Routing Problems: in this work, the half-way point is however adjusted after each exact pricing rounds, with similar settings.

\subsection{Design of PathWise}\label{Design}

\begin{figure}[ht]
	\centering
	\includegraphics[width=1\textwidth]{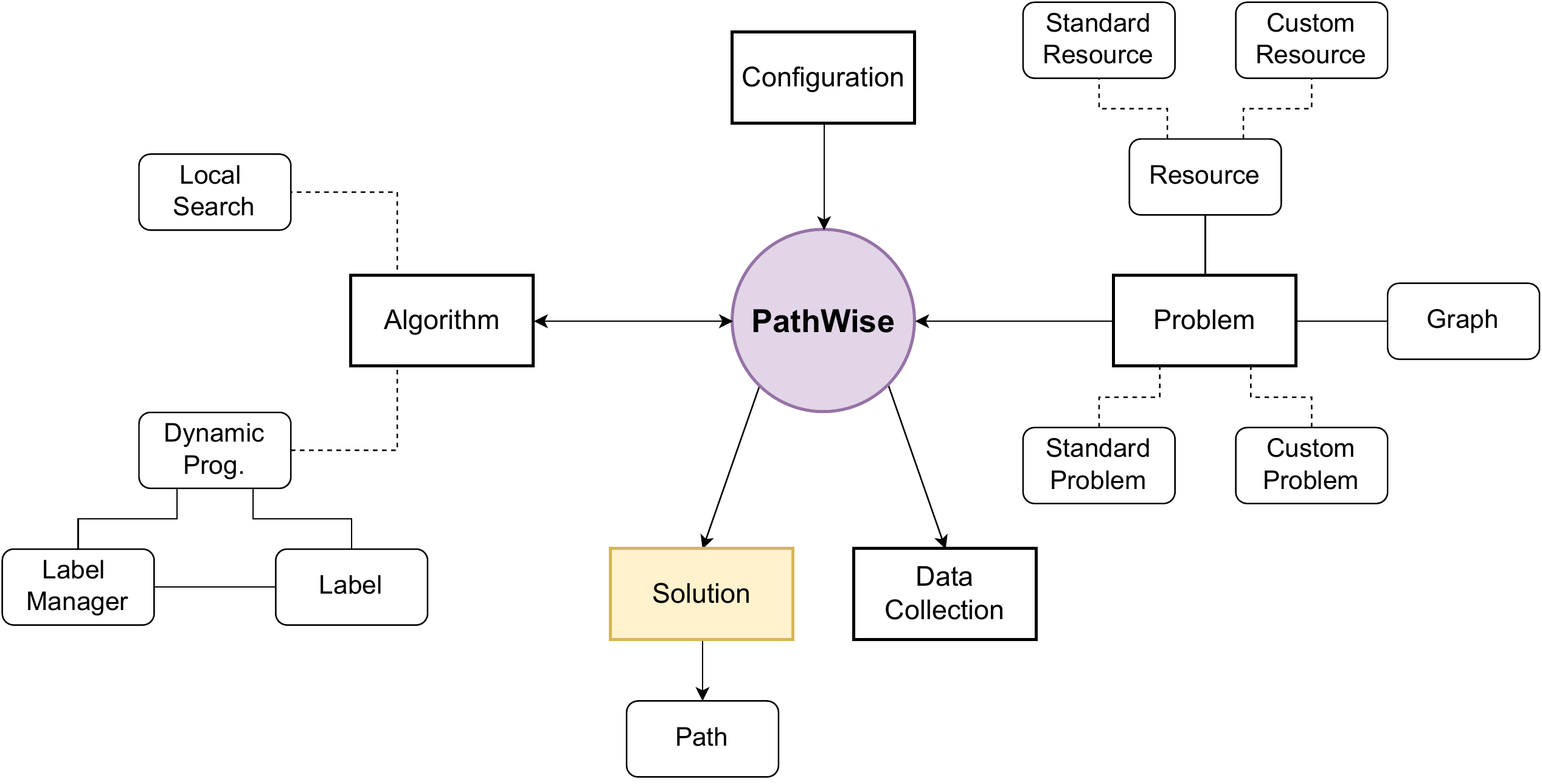}
	\caption{PathWise architecture}
	\label{figure:pw}
\end{figure}

In this section we propose the architecture of our library. We designed PathWise around a central solver unit that manages and interacts with other 5 major modules: configuration, problem, algorithm, solution and data collection. An example of the overall outline is sketched in Figure \ref{figure:pw}. At the beginning, PathWise communicates with the configuration component and gathers information to complete the initial setup of the library and the other modules. 
Problem data is presented to the framework by the user, either through instance files or directly in custom formulations, and is collected and arranged in the respective unit. 
Once a problem as been specified, PathWise then exploits one or more algorithms to solve, possibly at the same time, the input problem. Standard solutions are provided in the form of \emph{Paths} and additional statistics are collected and available both as output and possibly in the form of logs. 
\\We detail our main components in the following.

\paragraph*{Configuration} The Configuration module contains all the information that allows the framework to correctly complete setup, along with algorithm and output customization. For example, there are options to change the verbosity of the solver and data collection, choose aggressive or loose data management policies, activate parallel solving and select resolution algorithms. Furthermore, algorithms can still be additionally tuned by choosing relaxation techniques, extension and join strategies and selecting the critical resource in the bidirectional setup.
Default configurations are provided for both cyclic and acyclic networks, but they can be easily overridden by a user through a parameters file,  \emph{``pathwise.set''}, in the main folder.


\paragraph*{Problem} The Problem component stores and arranges in suitable data structures all the information about an instance, while implementing interfaces and efficient data access for the other units. We describe the problem with network information, an arbitrary number of resources and objective data.

A dedicated \emph{Graph} class is used to store and manage topological information about the network, such as node coordinates and edges, and can provide basic pre-processing operations and compute different adjacency lists to support dynamic programming algorithms. More in detail, we encode the presence of edges in incomplete networks with boolean variables, by either using a list of bitsets or a list of unordered maps, depending on the size of the instance. The former performs akin to an adjacency matrix: it allows direct access to information to guarantee maximum performance. The latter is used when facing very large, sparse graphs to save up memory, with look ups in constant time, on average.
Edge encoding is disabled by default for complete networks, but custom adjacency lists can still be provided, if needed. We remark that data structure configuration can be overridden by the user, through the parameters file.
The \emph{Resource} class allows to define and manage resource constraints by detailing the behaviour during extension and join procedures, for bidirectional algorithms, and feasibility checks. We designed a set of ready to use standard resource types from the literature and support custom, user outlined ones. In particular, we provide capacity and time constraints, limits for the number of nodes in a path and time windows.
A dedicated data structure, \emph{ResourceData} is used, for each resource, to collect lower and upper bounds, node and arc information: similarly to the graph topology, we implemented options to store arc consumption with either adjacency lists or unordered maps, depending on the size of the network. We use this data structure to also describe costs and profits for the objective function.
The first version of the library can solve problems with any number of monotonic resources, however, support for complex, non-monotonic ones is planned as well.

PathWise can handle both standard and custom problems. We define a problem as \emph{standard}, when all the resource types are among the ones supported by the library, and information is collected through a formatted instance file. In this case, graph and resource setup is performed automatically. 
We consider a problem as \emph{custom}, when either information requires ad-hoc retrieval methods, with possibly specialized data structures or pre-processing, or custom resources have to be defined. In these scenarios, PathWise provides simple interfaces to manually setup the network and resources. We report some examples in Section \ref{Usecase}. 

\paragraph*{Algorithm} The Algorithm module defines a suite of algorithms to solve the RCSPP, that may be heuristic or exact in nature. They are setup by the library, after configuration and problem definition have been completed and can be run either sequentially or in parallel.  For example, one may possibly run heuristic methods first, like a greedy or local search based algorithm, and then use the solution found as input to pre-processing techniques that might reduce the complexity of the problem, before running exact procedures.

This version of the library focuses on dynamic programming algorithms to find exact solutions: in particular, we propose a \emph{Bidirectional Dynamic Algorithm}, featuring the techniques defined in \ref{Algorithms}. However, the \emph{Algorithm} class presents interfaces that allow an expert user to describe dedicated custom methods of any kind. 

Our algorithms are designed to define the control flow, termination checks and, in general, all the steps required to obtain a solution without concerning about data structure implementation. They also take care of relaxation techniques and update them, and translate solutions in Paths. More in detail, we encoded states with a \emph{Label} data structure and defined a \emph{Label Manager} unit to handle them.

The \emph{Label} stores all the information regarding a partial path, such as the associated node, the direction (either forward or backward), the predecessor and the accumulated cost. Resource data is encoded as a vector of values, a \emph{snapshot} of all the consumption gathered in that particular state. Finally, we use bitsets to memorize both the nodes that have already been visited and the ones that are unreachable.

The \emph{Label Manager} unit handles labels, while performing core operations and providing encapsulation. That is, in our implementation, the algorithm does not have direct access to pools containing labels. Instead, the Label Manager answers requests from algorithms to find, return and possibly perform the extension of candidate labels towards non visited nodes, while inserting new ones in data structures and dealing with dominance checks. Additionally, it takes care of finding the best forward and backward labels that can be joined, when using bidirectional algorithms.

Our default implementation of the label manager presents multiple extension and join strategies, that can be selected through the parameters file. In particular, during extension steps, the algorithm can require all the labels of a certain node (Node selection) or only the top label from every node in a round robin setting (Round Robin selection), based on a specified metric. The join procedure can be performed in a naive fashion, for testing purposes, by trying to combine all the possible forward labels with all the backwards, or can exploit primal bounding to discard sub-optimal pairings.  

Overall, a template is present to allow the user to define specialized, custom version of dynamic algorithms too, while being able to fully exploit already defined data structures for labels and all the functionalities proposed by the Label Manager. Alternatively, even the Label Manager can be customized, while interfacing with either standard or custom algorithms.

Finally, we designed and implemented functionalities to support the user with algorithm debugging, allowing the search of specific labels in data structures and the generation of every label of a specified tour.


\paragraph*{Solution} The Solution component manages all the data structures required to handle output.  
In the standard scenario, a solution is provided by an algorithm and stored as a \emph{Path}, that is, a dedicated class which collects the objective cost, a tour of visited nodes, resource consumption and additional insights about its quality. Custom data structures for ad-hoc problems that can be modelled as RCSPP can be employed as well.
The solver can collect a set of Paths, depending on the number of algorithms used and their settings, rank and organize solutions over different metrics and provide a selection either as output for the user or as input for other algorithms. 

\paragraph*{Data Collection} Additional information is saved by the \emph{Data Collection} unit. This module can fetch features regarding the configuration, algorithm runs and the obtained Paths to show reports and store logs. For example, global and routine optimization times can be profiled, along with statistics such as the number of algorithm and relaxation iterations, insertion and join data, and more broadly, label generation details. We designed this unit with generic methods, that can be included in custom algorithms, label managers and problems as well, to profile their behaviour. 

We envision three main uses of the collected data. First, gaining insights on the solution and the behavior of the algorithms over the problems users are solving.
Second, algorithm profiling to analyze standard and custom modules of the library.
Third, supporting the development of possibly new data driven approaches. However, we note that data collection can be completely turned off in the parameters, if needed.

\subsection{Use cases of PathWise} \label{Usecase}
In the following, we summarize and provide short examples of the main use cases of PathWise.

\paragraph*{Standard Problem} The basic use of PathWise is to solve standard RCSPP problems stored in instance files of compatible format. PathWise can be used as a standalone executable, loading the instance and writing the output to a file or integrated in the user's application with few lines of code. The problem can be composed of an arbitrary set of resources of standard type. The pseudo-code reported in \nameref{ref:sp} provides an example of the integration of PathWise in the user's application. 

\begin{lstlisting}[language=C++, caption={Standard Problem}, title={Standard Problem}, label={ref:sp}]
    //Create a PathWise instance
    PathWise pw = PathWise(<file_name>);
    //Solve the problem
    pw.solve();
    //Get solutions 
    pw.getSolutions();
\end{lstlisting}

\paragraph*{Custom problem with standard resources} 
PathWise can be used to solve non-standard problems when, for example, data is organized in files of non compatible format or when resource consumption data needs dedicated data management.

The user is required to extend the base class of a problem and provide the data management components. The pseudo-code reported in \nameref{ref:cp2} provides an example of the definition of an ad-hoc problem extending the Problem base class. The user needs to override the method reading the network from file and manage the creation of the related resources implementing the method \textit{readProblem} accordingly.

\begin{lstlisting}[language=C++, caption={Custom Problem}, title={Custom Problem}, label={ref:cp2}]
//Create a Custom Problem
class CustomProblem: public Problem {
public:
    CustomProblem() = default;
    ~CustomProblem() = default;

    // Overridden methods
    void readProblem(std::string file_name);
};
\end{lstlisting}

\paragraph*{Custom problem with custom resources} 
The user may want to add a custom resource with dedicated extension and feasibility check methods. The user needs to extend the Resource base class and adapt the data management component of the problem class as described above. 
The pseudo-code reported in \nameref{ref:cr} provides an example of the definition of a  resource with custom consumption, e.g. non-linear. For brevity, in the pseudo-code \textit{c\_n, c\_v, c\_v\_fw, c\_v\_bw, dir} stand for current node, current value, current value forward, current value backward and direction, respectively.
The user is required to set the proper bounds in the \textit{init} method, to provide updating mechanisms in the \textit{extend} method, to determine when the current value at a given node of the network is feasible in the \textit{isFeasible} method. The \textit{join} method is specifically conceived for bi-directional algorithms when the current value is obtained as the union of two values in opposite directions.

\begin{lstlisting}[language=C++, caption={Custom Resource}, title={Custom Resource}, label={ref:cr}]
//Create a custom non-linear resource
class NonLinearCapacity: public Resource {
public:
    NonLinearCapacity() = default;
    ~NonLinearCapacity() = default;

    // Overridden methods
    void init(int origin, int destination);
    float extend(float c_v, int i, int j, bool dir);
    bool isFeasible(float c_v, int c_n, float bounding, bool dir);
    float join(float c_v_fw, float c_v_bw, int i, int j);
};
\end{lstlisting}

\paragraph*{Custom algorithm} The more advanced user may be interested in providing a custom algorithm to solve either a standard problem or a custom problem. The user needs to extend the \textit{Algorithm} base class providing the implementation of the main methods. The most relevant components of the abstract class are reported in pseudo-code \nameref{ref:ca}. In particular the \textit{solve} method is the main entry point of the algorithm and the method \textit{getSolutions} is invoked to retrieve the pool of computed solutions.

\begin{lstlisting}[language=C++, caption={Custom Algorithm}, title={Custom Algorithm}, label={ref:ca}]
class CustomAlgorithm: public Algorithm {
public:
    // Overridden methods
    // Basic data query-set methods not reported here
    ...
    inline int getStatus() {return status;}

    // Solve methods
    virtual void solve(Problem* problem) = 0;
    virtual std::vector<Path> & getSolutions() = 0;
};
\end{lstlisting}

\paragraph*{Custom dynamic programming algorithm} Advanced users that want to customize some components of the dynamic programming algorithm do not need to write a dedicated algorithm and extend the \textit{Algorithm} base class.  
They may extend the \textit{DpAlgo} base class. The most relevant components of the abstract class are reported in pseudo-code \nameref{ref:dp}. The user needs to implement the main methods \textit{solve} and \textit{getSolutions} as before but can leverage on the already implemented \textit{LabelManager} component.

\begin{lstlisting}[language=C++, caption={Custom DP Algo}, title={Custom DP Algo}, label={ref:dp}]
class CustomDPAlgo: public DPAlgo {
public:
    // Overridden methods
    // Basic data query-set methods not reported here
    ...
    inline int getStatus() {return status;}
    
    // Label manager
    inline void setLabelManager(LabelManager* lbl_manager){this->lbl_manager = lbl_manager;}
    inline void setLabelManager(){lbl_manager = new LabelManager;}
    
    // Solve methods
    ...
};
\end{lstlisting}

\paragraph*{Custom label manager} An interesting advanced customization of the dynamic programming algorithm is to provide a different label manager component. Indeed, the way labels are generated, stored, maintained and deleted is at the heart of the dynamic programming algorithm.
The advanced users may want to experiment new ways to perform these basic operations on labels.
To do so, users need to extend the \textit{LabelManager} base class.  
The most relevant components of the class are reported in pseudo-code \nameref{ref:lm}. 

\begin{lstlisting}[language=C++, caption={Custom Label Manager}, title={Custom Label Manager}, label={ref:lm}]
class CustomLabelManager: public LabelManager{
public:

    // Overridden methods
    // Startup
    void init(Problem* problem);

    // Get candidate
    bool candidatesAvailable(bool direction);
    Label* getCandidate(bool direction);

    // Label extension
    bool isExtensionFeasible(Label& current_label, int id);
    void extendLabel(Label* current_label, Label & new_label, int id);

    // Label insertion
    Label* insert(Label& new_label);

    // Join
    void join();

    // Get Solution
    PW_LabelPair getSolutionLabels();
};
\end{lstlisting}

The use cases provided in this section represent some short meaningful examples of different levels of user's interaction with PathWise. More comprehensive and commented tutorials are shipped with the library.

\section{Computational Experiments}

We report, in this section, the computational analysis of the framework, designed around 3 different class of problems, with the aim of evaluating the performance and flexibility of PathWise when facing different scenarios.  We considered both cylic and acyclic networks, with single and multiple resource constraints.

\subsection{Test-bed instances}

Experiments were performed on 3 datasets, including both instances found in the literature and newly generated ones. 

\paragraph*{SPPRCLIB}
The first dataset was taken from the online repository \citep{spprclibrepo}. It consists of 45 instances, describing complete graphs, with up to 262 nodes, with a single capacity constraint. They were derived from CVRP instances that have been solved through column generation \citep{Fukasawa2006, Jepsen2008}, and present both positive distance costs on arcs and positive prizes (i.e., negative costs) on nodes, thus making the problem cyclic. They are categorized with labels (A, B, E, G, M, P) according to the original authors. We used this set as reference for performance on shortest path problems with a single resource constraint and cyclic networks.

\paragraph*{DIMACS}
The second dataset was selected instead from the 9th DIMACS Challenge \citep{DIMACS9th}, that collects several USA road networks. In this case, instances are characterized by acyclic, large (up to several millions of nodes) sparse graphs, with a single time constraint. 
For this set, we selected, overall, 160 instances from New York (NY), San Francisco Bay Area (BAY), Colorado (COL) and California and Nevada (CAL) networks, with the same destinations and upper bounds that are reported in \citep{Lozano2020}. Other instances were also considered for preliminary experiments. We used this dataset to profile PathWise performance on shortest path problems with a single resource constraint and acyclic graphs.

\paragraph*{Prize Collecting (PC)} 
Finally, for the third dataset, we generated a new problem starting from CVRPLIB instances \citep{CVRPLIB}. Specifically, we designed an application of an RCESPP on a complete, cyclic graph (up to 100 nodes), in which every arc of the network presents a negative cost. Therefore, the optimal solution maximizes the collection of prizes within the feasible resource consumption. Similar problems in the literature can be found in the form of variants of arc orienteering problems \citep{GAVALAS} and prize collecting TSP \citep{Balas}.
However, in our setting, multiple resource consumptions are present and the problem is constrained by two capacity bounds, a node limit and time windows.
In fact, we generated 48 new instances that differ from each other in the number of nodes, the capacity thresholds and the problem node limit. 
We designed this dataset to test PathWise performance and modelling abilities when facing shortest path problems with multiple resource constraints.

Additional details regarding instance generation are reported in the Appendix, in Section \ref{instgen}.

\subsection{Configuration}
We developed PathWise in C++17. Although PathWise was designed as a standalone package, this pre-release implementation of the library exploited Boost 1.74 to manage dynamic bitsets. Compilation was performed through GCC 11.2 with the ``O3'' optimization flag and tests were executed on a machine running Kubuntu 22.04 that was equipped with an eight-core Intel i9-11900 @ 2.50 GHz and 32GB RAM. 

After preliminary experiments, we found it preferable to use two distinct set of settings for cyclic networks (RCESPP) and acyclic ones (RCSPP). In the former case, we exploited bidirectional parallel labeling, with an initial even split for the critical resource, and a node first candidate selection approach during extensions. More in detail, we chose the node presenting the minimum cost label in the pool. 
We also considered the four different relaxation techniques that guarantee an elementary solution, reported in Section \ref{Algorithms}:
\begin{itemize}
\item {\tt DSSR} 
\item {\tt DSSRC} 
\item {\tt NG-DSSRC}
\item {\tt NGC-DSSRC}
\end{itemize}

Experimentally, we found out that using an initial neighbourhood of size 16, for ng-route relaxations, was a good compromise between solution quality and performance, for this set of tests.

In the latter case (RCSPP), we used, sequentially, bidirectional labeling, to better control primal bounding, and a round robin candidate selection strategy for extension. In this scenario, we selected and extended, in turn, the minimum cost label of each node. Additional settings for memory management were implemented by disregarding label information related to already visited nodes and by storing sparse data through unordered maps. \\
Finally, we take advantage, in both setups, of bounded join procedures.


\subsection{Experiments}
In the following, we report our main results. When available, we also present solutions from other works in the literature that tackle the same problems and instances.

\paragraph*{Single resource constraint, cyclic network profiling.}
In our first round of experiments, we profile the performance of PathWise when solving the RCESPP on complete, cyclic graphs, with a single resource constraint.
In Table \ref{tab:spprclib} we present the absolute time in seconds when solving every instance of SPPRCLIB with PathWise while exploiting different relaxations. Additionally, we also report results when using Branch-and-Cut (B\&C), from \citep{spprclib}. We remark, however, that PathWise and B\&C experiments run on different machines, and therefore the comparison serves only as a reference of the general performance of specialized techniques.
 The best times are marked in bold.  Experiments with missing results hit the 1 hour time limit.
\begin{table*}[p]
    \small
    \begin{tabularx}{\textwidth}{X*{6}{r}}\toprule
        & & \multicolumn{4}{c}{PathWise}\\
        \cmidrule(lr){3-6}
        Instance	&B\&C	& DSSR	& DSSRC	& NG-DSSRC	& NGC-DSSRC\\
	\midrule
        A-n54-k7-149	&6.96	&\textbf{0.66}	&0.77	&0.83	&0.70\\
        A-n60-k9-57	&36.55	&\textbf{1.90}	&2.00	&1.91	&1.94\\
        A-n61-k9-80	&4.44	&\textbf{0.49}	&0.66	&1.26	&0.59\\
        A-n62-k8-99	&17.94	&\textbf{0.30}	&0.55	&2.48	&0.51\\
        A-n63-k9-157	&3.16	&\textbf{0.17}	&0.28	&1.30	&0.28\\
        A-n63-k10-44	&2.12	&\textbf{0.13}	&0.32	&0.35	&0.29\\
        A-n64-k9-45	&14.57	&\textbf{0.58}	&1.51	&5.37	&1.92\\
        A-n65-k9-10	&4.43	&\textbf{0.16}	&0.29	&1.04	&0.29\\
        A-n69-k9-42	&1.76	&\textbf{0.22}	&0.34	&0.44	&0.33\\
        A-n80-k10-14	&12.14	&\textbf{2.25}	&4.67	&15.21	&5.11\\
	\addlinespace
        B-n45-k6-54	&\textbf{1.32}	&5.35	&1.88	&6.21	&1.75\\
        B-n50-k8-40	&11.01	&\textbf{0.17}	&0.19	&1.76	&0.19\\
        B-n52-k7-15	&\textbf{1.00}	&1.68	&1.14	&11.83	&1.39\\
        B-n57-k7-20	&\textbf{1.74}	&-	&-	&-	&-\\
        B-n66-k9-50	&66.93	&\textbf{0.32}	&1.45	&20.35	&1.42\\
        B-n67-k10-26	&4.62	&\textbf{0.18}	&0.47	&0.80	&0.47\\
        B-n68-k9-65	&11.88	&\textbf{0.59}	&2.29	&12.66	&2.10\\
        B-n78-k10-70	&24.30	&\textbf{1.13}	&2.31	&6.31	&2.68\\
	\addlinespace
        E-n76-k7-44	&6.02	&\textbf{0.45}	&0.62	&6.61	&0.69\\
        E-n76-k10-72	&1.19	&\textbf{0.43}	&0.50	&3.00	&0.53\\
        E-n76-k14-102	&14.77	&0.87	&1.27	&\textbf{0.23}	&1.27\\
        E-n76-k15-40	&19.59	&0.65	&1.04	&\textbf{0.19}	&1.04\\
        E-n101-k8-291	&8.08	&\textbf{0.85}	&1.04	&2.07	&0.97\\
        E-n101-k14-158	&37.84	&\textbf{1.42}	&1.73	&1.72	&1.74\\
	\addlinespace
        G-n262-k25-316	&\textbf{53.00}	&-	&-	&-	&-\\
	\addlinespace
        M-n101-k10-97	&3.12	&\textbf{2.00}	&7.23	&83.91	&6.08\\
        M-n121-k7-260	&\textbf{34.46}	&-	&-	&-	&-	\\
        M-n151-k12-15	&\textbf{78.03}	&523.95	&498.12	&2109.75	&480.95\\
        M-n200-k16-143	&\textbf{3.18}	&2089.31	&779.33	&2841.24	&752.01	\\
        M-n200-k17-12	&\textbf{17.75}	&207.88	&482.52	&2295.81	&499.98	\\
	\addlinespace
        P-n50-k7-92	&2.42	&\textbf{0.15}	&0.31	&0.21	&0.31\\
        P-n50-k8-19	&0.36	&\textbf{0.28}	&0.49	&1.06	&0.48\\
        P-n50-k10-24	&0.72	&\textbf{0.04}	&0.05	&\textbf{0.04}	&0.05\\
        P-n51-k10-30	&2.18	&0.17	&0.27	&\textbf{0.06}	&0.27\\
        P-n55-k7-116	&0.58	&\textbf{0.04}	&0.06	&0.56	&0.07\\
        P-n55-k8-260	&1.20	&0.18	&\textbf{0.16}	&0.32	&\textbf{0.16}\\
        P-n55-k10-44	&2.14	&0.15	&0.22	&\textbf{0.09}	&0.22\\
        P-n55-k15-88	&3.97	&0.19	&0.27	&\textbf{0.04}	&0.30\\
        P-n60-k10-24	&1.04	&\textbf{0.05}	&0.11	&0.29	&0.12\\
        P-n60-k15-8	&1.95	&0.11	&0.16	&\textbf{0.02}	&0.15\\
        P-n65-k10-102	&6.65	&0.59	&0.92	&\textbf{0.43}	&0.90\\
        P-n70-k10-12	&0.24	&\textbf{0.11}	&0.22	&0.85	&0.21\\
        P-n76-k4-41	&\textbf{1.85}	&27.11	&73.55	&272.12	&60.51\\
        P-n76-k5-16	&\textbf{0.57}	&51.11	&14.10	&18.29	&15.38\\
        P-n101-k4-174	&\textbf{11.25}	&104.32	&51.03	&1659.42	&37.30\\
	\hline
	\addlinespace
	Average	&\textbf{10.75}	&71.33	&45.99	&188.51	&45.03\\
	\bottomrule
    \end{tabularx}
    \caption[Comparison on SPPRCLIB instances]{Comparison on SPPRCLIB instances. We report the absolute time [s] for a state of the art implementation of Brach-and-Cut and PathWise, when using different relaxations.} \label{tab:spprclib}
\end{table*}

We can observe that PathWise is competitive on most instances, scoring the best results in most of them. Although the comparison with Branch-and-Cut is only sketched, our framework includes a dynamic programming algorithm that generally performs in the same order of magnitude or better, struggling only in specific scenarios. As expected, it vastly outperforms dynamic programming implementations that do not exploit relaxation techniques: for example, the one presented in \citep{spprclib} hits almost systematically timeout in all experiments.  However, some intrinsic weaknesses remain when facing instances with a high number of nodes or when too many arcs with negative cost are present in the network. Indeed, in these scenarios the number of generated labels cannot be easily controlled. This is expected and in line with the literature \citep{spprclib}: improving efficiency remains an open research direction. \\
PathWise different configurations can be instead directly compared. Results here are mixed since {\tt DSSR} obtains the best results in about 73\% of the solved instances whilst {\tt NGDSSR} in about 22\% of the cases. The other two configurations seldom produce the best score but perform on average better than {\tt DSSR} and {\tt NGDSSR}, thus indicating that they might be able to better manage labels, in particular when facing large networks or difficult instances. 
Overall, these results suggest that different relaxations might be best suited to different instances. In fact, in these experiments alone, choosing the best performing relaxation instead of the worst one would provide, on average, about 73\% faster solutions: being able to accurately predict which technique should be used for a given problem would have important impact on performance. 

\paragraph*{Single resource constraint, acyclic network profiling. }In our second round of tests, we study how PathWise behaves on the large, sparse, acyclic networks of the DIMACS dataset, when solving the RCSPP with a single resource constraint. In Table \ref{tab:dimacs}, we report the average time (Time) in seconds and the overall number of solved instances for PathWise and for the state-of-the-art A* based algorithm RC-EBBA* \citep{EBBA}. In this case, we compiled and tested the algorithm on our machine, using the corresponding repository available online \citep{A*REPO}. Best results are marked in bold. Timeout was set to 4 hours.

\begin{table*}[h!]

    \small
    \begin{center}

    \begin{tabularx}{0.85\textwidth}{X*{6}{r}}\toprule
        & \multicolumn{2}{c}{PathWise} & &\multicolumn{2}{c}{RC-EBBA*} \\
        \cmidrule(lr){2-3} \cmidrule(lr){5-6} 
        Network	&Time	&Solved	&&Time	&Solved	\\
	\midrule
        NY	&69.58	&40/40&	&\textbf{0.02}	&40/40	\\
        BAY	&84.36	&40/40&	&\textbf{0.02}	&40/40\\
        COL	&311.08	&40/40&	&\textbf{0.04}	&40/40	\\
        CAL	&1232.93	&38/40&	&\textbf{6.13}	&\textbf{40}/40	\\
	\hline
	\addlinespace
        Overall& & 158/160 && & \textbf{160}/160 \\ 
	\bottomrule
    \end{tabularx}
    \end{center}
    \caption[Comparison on DIMACS instances]{Comparison on DIMACS instances. We report average time [s] and number of solved instances for each network type and algorithm.} \label{tab:dimacs}
\end{table*}

Performance on SPPRC instances with a single resource constraint is not competitive with the latest algorithms from the literature. While PathWise never hits the time limit and solves all small sized instances but two without incurring in memory issues,  solution time is not comparable to specialized algorithms such as RC-EBBA*.
At this time, more ad-hoc optimizations are required to tackle this specific class of problems, while providing similar speed-ups and good memory management. To this end, an implementation of an exact A* algorithm in our framework would be beneficial to performance and help when facing larger instances.
 However, we remark that our library, being generic in nature, can be used to deal with problems with any amount and type of resource constraints. We expect advanced A* based techniques to possibly provide a speed-up in these settings as well, however, the impact on run time for problems that require more complex dominance checks needs to be verified. We also note that these techniques might be unsuitable for some classes of problems: for example, instances presenting only non-monotonic custom resources might not be solvable in an efficient way or would at the very least require a tailored, much more involved pre-processing. 


\paragraph*{Multiple resource constraints.}
In our third round of experiments, we profile how PathWise performs on RCESPP when solving PC instances, featuring complete networks and multiple resource constraints. Indeed, we were able to easily model this setting with our library, since any problem would be automatically configured, as long as resource constraints were among the defined ones and data was correctly presented.  This also allowed to turn on and off particular resource constraints when needed: we used this feature to evaluate different setups and to choose the experimental configuration reported in the following.
In Table \ref{tab:lp}, we present the average time (Time) in seconds and the overall number of unsolved instances before timeout for PathWise, when using different configurations. We detail results for each each instance size (n), capacity threshold (C) and node limit (NL). Best scores are marked in bold. For these experiments, timeout was set to one hour.

\begin{table*}[h!]
    \footnotesize
    \begin{tabularx}{\textwidth}{X*{11}{r}}\toprule
        & & & \multicolumn{2}{c}{DSSR}& \multicolumn{2}{c}{DSSRC} & \multicolumn{2}{c}{NG-DSSRC}&\multicolumn{2}{c}{NGC-DSSRC}\\
        \cmidrule(lr){4-5} \cmidrule(lr){6-7} \cmidrule(lr){8-9} \cmidrule(lr){10-11}  
        n	&C	&NL	&Time &Uns. &Time &Uns. &Time &Uns. &Time &Uns.\\
        \midrule
        \multirow{5}{*}{50}	&\multirow{2}{*}{25} 
            &8	&0.32	&0	&0.35	&0	&\textbf{0.24}	&0	&0.35	&0\\
            & &18	&\textbf{0.41}	&0	&1.01	&0	&0.50	&0	&1.12	&0\\
        \addlinespace
	                    &\multirow{2}{*}{40}	
            &8	&9.28	&0	&12.37	&0	&\textbf{8.16}	&0	&12.42	&0\\
            &&18	&\textbf{14.65}	&0	&39.10	&0	&26.37	&0	&42.99	&0\\

                                      \midrule
        \multirow{5}{*}{100}	&\multirow{2}{*}{25}	

        &8	&4.82	&0	&6.96	&0	&\textbf{2.35}	&0	&7.27	&0\\
        &&18	&\textbf{3.89}	&0	&6.81	&0	&4.83	&0	&7.38	&0\\
                                                \addlinespace
	                           &\multirow{2}{*}{40}	
        &8	&276.28	&0	&489.68	&0	&\textbf{144.12}	&0	&493.50	&0\\
	&	&18	&1147.94	&\textbf{0}	&173.49	&3	&104.55	&3	&179.76	&3\\
          
                                                	\hline
	\addlinespace
Overall		
&	&	&182.20	&0	&85.74	&3	&31.85	&3	&87.32	&3\\
\\

	\bottomrule
    \end{tabularx}
    \caption[Comparison on PC instances, with TW]{Comparison of PathWise on PC instances, when using different relaxations. We report the average time [s] and the number of timeouts for each instance class.} \label{tab:lp}
\end{table*}

Overall, the amount of negative arcs makes this problem challenging for labeling algorithms, since finding an elementary path is inherently much harder. 
Nonetheless, results show that all configurations can solve most instances without hitting the time limit. Performance scales as expected: when the problem is well constrained, PathWise can solve instances within seconds. However, when the number of nodes is large and capacity thresholds are loose, the problem becomes much more hard to tackle. This is in line with the literature and, under these specific conditions, some timeouts occur for all configurations but {\tt DSSR}. More in detail, {\tt DSSR} seems to perform better than other relaxations when the node limit is higher, whilst {\tt NG-DSSRC} is the better algorithm on shorter paths. The other two configurations are instead dominated. Again, this is an another indicator that different relaxation types seem to be better suited to particular classes of instances.

\paragraph*{Extension strategies profiling.}
Finally, in Table \ref{tab:selection} we study the impact of exploiting different label selection strategies during the extension step of the labeling algorithm. In particular, we present PathWise results when solving SPPRCLIB, DIMACS and PC instances. We consider two candidate selection strategies: Round Robin selection and Node selection. 
 For each problem and selection policy, we present the average time (Time) for solved instances by both configurations and the overall number of unsolved instances (Unsolved). Additionally, we report if the network is cyclic or acyclic (Network). Best results are marked in bold. We kept the same timelimit used in previous experiments and used the {\tt DSSRC} configuration to tackle cyclic networks.

\begin{table*}[h!]
\small
\begin{center}
    \begin{tabularx}{0.9\textwidth}{X*{8}{r}}\toprule
        &&& \multicolumn{2}{c}{Round Robin sel.} & &\multicolumn{2}{c}{Node sel.}\\
        \cmidrule(lr){4-5} \cmidrule(lr){7-8}
        Problem	& Network &&Time	&Unsolved &	&Time	&Unsolved\\
	\midrule
        SPPRCLIB &Cyclic&		&82.55	&4	&	&\textbf{28.22}	&\textbf{3}\\
        \addlinespace
        DIMACS  &Acyclic&       	&414.25	&\textbf{2}	&	&396.70	&17\\
        \addlinespace
        PC 	&Cyclic &           &190.75
	&3	&	&\textbf{85.74}	&3\\
\bottomrule
\end{tabularx}
\end{center}
    \caption[Extension strategies profiling]{PathWise profiling over SPPRCLIB, DIMACS and PC instances when using different extension strategies. We report the average time [s] and the number of unsolved instances.} \label{tab:selection}
\end{table*}

Overall, node selection was, on average, the best configuration when facing problems with cyclic networks. In this setting, solution time was lower and the number of unsolved instances was smaller. However, when facing large acyclic networks this strategy hit multiple timeouts, whilst Round Robin selection performed better. Experimentally, extending all the available labels for a particular node in an acyclic graphs makes label generation much harder to control in an efficient way. Indeed, being able to automatically configure the algorithm behavior depending on the instance class or its features seems promising, and has been studied in the literature \citep{Schede_2022} in other settings, and would allow to obtain further performance improvements. We expect this hold true in many other scenarios. For example, NG-route based relaxation techniques could be additionally tuned by selecting a custom neighbourhood size according to instance properties. 

\section{Concluding remarks}

In this paper, we described the design and implementation of \textit{PathWise}, a flexible, open-source library for the solution of the Resource Constrained Shortest Path Problem. 

PathWise is primarily an easy to use library for the standard user, shipped with state of the art algorithms. 
Furthermore, the user can easily model multiple resources with little effort. Standard resources provided off-the-shelf have the potential to cover most of the applications, while new standard resources are expected to be added with subsequent releases of the library together with explanatory examples and tutorials.
PathWise is flexible for more advanced users' customizations allowing for ad-hoc representation of non-standard problems.
For both standard and custom problems, the user is provided with flexible data collection methods that permit to profile the solution and the behavior of the algorithms.

More in detail, the first release of PathWise, outlined in this paper, is shipped with an exact bi-directional dynamic programming algorithm implementing state-of-the-art techniques such as multiple relaxations, dynamic half-way point, memory compression techniques, different variants of label extension and join mechanisms, while presenting an easy to use configuration interface.
We tested PathWise on three representative classes of RCSPP. Namely, RCSPP on cyclic networks (SPPRCLIB), RCSPP on large acyclic networks (DIMACS) and RCSPP on ad-hoc cyclic networks (PC).
PathWise results competitive on the problems with cyclic networks instances, whilst it does not compare favorably against ad-hoc, less flexible algorithms when facing acyclic, large networks.

Computational results show that there is no dominant algorithm that suits all classes of problems. Instead, careful algorithm selection and configuration based on classes of instances or their features can have a positive impact on solution times. 
This is a clear indication that data driven methods could possibly help in finding and setup the best algorithm, therefore improving performance. 

We have a clear vision for the future releases of PathWise. First, we want to fill the performance gaps with ad-hoc algorithms for large acyclic networks while maintaining the flexible philosophy of the library and improving the developed algorithms with high level parallelization strategies.
Second, we intend to further expand the library with the implementation of complex, non monotonic resources, eventually using them in a similar fashion to standard ones. Finally, we intend to add data driven methodologies for automatic algorithm selection and configuration.

\bibliography{bibliography.bib}
\newpage
\appendix
\appendix
\section{Appendix}

\subsection{Prize Collecting instance generation} \label{instgen} 
In this section, we report the main characteristics regarding PC instance generation. Starting from the ``Loggi'' set (published on 2021) of CVRPLIB \citep{CVRPLIB}, we generated 48 new instances that differ in the number of nodes, the capacity thresholds and the problem node limit. 
More in detail, the original problem presents positive arc costs and a capacity constraint. After preliminary experiments, we found the following setup to produce meaningful experiments.
\begin{itemize}
\item For each instance, we selected either the first 50 or 100 nodes ($n$). We kept the network complete, and we reported the original arc costs, based on distances, as prizes, that is, with negative values.
\item  We kept the original node consumption for capacity but updated the upper bound to either 25 or 40 ($C$).
\item We added a second capacity constraint. We generated node consumption between 1 and 10, like the original values in CVRPLIB instances. Then, we set an upper bound between 80\% and 120\% of $C$.
\item We enforced a node limit constraint that restricts the maximum length of a path to either 8 or 18 nodes ($NL$).
\item Finally, we added a Time Windows (TW) constraint. We obtained times for arc traversal by dividing distances by 100, whilst service times for nodes were instead generated from a set of values (10, 20, 30 and 40).
We then enforced about 80\% wide time windows, to guarantee feasibility and non trivial instances, and 20\% narrow ones. Arrival time ($at$) was generated between 0 and 1000. Departure time ($dt$) as follows:

\begin{equation}{
    dt = 
    \begin{cases}
        at + 100*random(1,4) &\text{if wide TW} \\
        at + 100*random(0.1, 0.6) &\text{if narrow TW}\\
    \end{cases}
    }
\end{equation}

We also made sure that when arriving at a node there was always enough time to serve it.

\end{itemize}

We remark that generation was always performed randomly through a uniform distribution.

\end{document}